\documentclass[a4paper]{article}
\usepackage{pdfsync}
\usepackage[utf8]{inputenc}
\usepackage{graphicx}
\usepackage{amssymb,amsfonts,amsmath,amsthm}
\usepackage{hyperref}
\usepackage{multirow}
\usepackage{verbatim}
\usepackage[rgb,dvipsnames]{xcolor}
\usepackage{booktabs}
\usepackage[sectionbib,round]{natbib}
\usepackage{textcomp}
\usepackage{upquote}
\usepackage{bm}
\usepackage{dirtree}
\usepackage{pdflscape}

\theoremstyle{plain}

\theoremstyle{definition}

\begin{document}

\title{\textbf{Interdependence between Green Financial Instruments and Major Conventional Assets: A Wavelet-Based Network Analysis}}


\author{R. Ferrer$^1$, R. Ben\'{\i}tez$^2$, V. J. Bol\'os$^2$ \\ \\
{\small $^1$ Dpto. Econom\'{\i}a Financiera y Actuarial, Facultad de Econom\'{\i}a.} \\
{\small $^2$ Dpto. Matem\'aticas para la Econom\'{\i}a y la Empresa, Facultad de Econom\'{\i}a.} \\
{\small Universidad de Valencia. Avda. Tarongers s/n, 46022 Valencia, Spain.} \\
{\small e-mail\textup{: \texttt{roman.ferrer@uv.es}, \texttt{rabesua@uv.es}, \texttt{vicente.bolos@uv.es}}} \\}

\date{April 2021}

\maketitle

\begin{abstract}
This paper examines the interdependence between green financial instruments, represented by green bonds and green stocks, and a set of major conventional assets, such as Treasury, investment-grade and high-yield corporate bonds, general stocks, crude oil, and gold. To that end, a novel wavelet-based network approach that allows for assessing the degree of interconnection between green financial products and traditional asset classes across different investment horizons is applied. The~empirical results show that green bonds are tightly linked to Treasury and investment-grade corporate bonds, while green stocks are strongly tied to general stocks, regardless of the specific time period and investment horizon considered. However, despite their common climate-friendly nature, there is no a remarkable association between green bonds and green stocks. This means that these green investments constitute basically two independent asset classes, with a distinct risk-return profile and aimed at a different type of investor. Furthermore, green financial products have a weak connection with high-yield corporate bonds and crude oil. These findings can have important implications for investors and policy makers in terms of investment decision, hedging strategies, and sustainability and energy policies.
\end{abstract}

\section{Introduction}

Climate change has become one of the priority challenges that humanity faces in the twenty-first century. In order to avoid the devastating effects of climate change on human life and ecosystems, it is vital to adopt rapid and large-scale actions that reduce fossil fuels consumption and carbon emissions. 
Nowadays, green bonds and green stocks constitute the two major environmental-friendly financial instruments and are expected to play a critical role in mobilizing the staggering amount of capital required to finance the massive transformational projects needed to the transition to a low carbon economy.

Green bonds are a type of fixed income security whose essential difference from ordinary bonds is that the proceeds of green bonds are exclusively employed to finance environmental or climate-friendly projects in areas such as renewable energy, energy efficiency, clean transport, sustainable water management and green buildings. Since the first ever green bond was issued in 2008 by the European Investment Bank (EIB), the green bond market has become one of the fastest-growing segments of international capital markets, especially in the wake of the publication of the Green Bond Principles (GBP)\footnote{The Green Bond Principles (GBP) are a set of internationally recognized voluntary guidelines for the issuance of green bonds that recommend transparency and disclosure and promote integrity in the green bond market. The GBP have remarkably enhanced the level of transparency and credibility of the green fixed income market and also avoided green washing, acting as a primary catalyst for the extraordinary growth of this market over the past few years.} by the International Capital Market Association (ICMA) in January 2014.   In this regard, the global green bond issuance in 2014 was USD $36.6$ billion, while, according to Climate Bonds Initiative (CBI), issuance of green bonds in 2020 has reached USD $269.5$ billion (just above USD $266.5$ billion in 2019). Despite the dramatic slowdown in issuance during the first half of 2020 caused by the COVID-19 outbreak, the 2020 figure is the highest since the market's inception and consolidates a trend of nine consecutive years of growth of the green bond market. In addition, the cumulative green bond issuance since their first appearance in 2007 crossed over the USD $1$ trillion milestone in September 2020. However, despite this recent boom, the global green bond market still remains a very small fraction of the more than USD $100$ trillion overall fixed income market, which can be interpreted as a powerful signal of the enormous potential for growth of this green market.

Green stocks can be defined as shares of companies whose primary business has a beneficial impact on the environment. These companies are mainly concentrated in sectors such as alternative energy, energy and material efficiency, clean transportation, water management, and waste management and recycling. The green equity market has also witnessed an unprecedented surge since the 2015 Paris Climate Agreement due to the massive capital flows entering into green stocks in a context where investors seek environmentally responsible investments that have at the same time a high potential of return. More specifically, according to Morningstar, investors poured USD $51.1$ billion of new money into US environmental, social and governance (ESG) investment funds in 2020, more than doubling the inflows of USD $21.4$ billion in 2019. This figure represents the fifth consecutive annual record and is particularly remarkable taking into consideration the sharp market downturn fueled by the Coronavirus during the first half of 2020. Morningstar data also reveal that ESG funds made up about a fourth of all the money that flowed into US equity and bond mutual funds in 2020. As an example, ESG funds only accounted for about 1\% of inflows into total US funds in 2014.

The aim of this paper is to examine the degree of interconnection between green financial instruments, which are represented by green bonds and green stocks, and a set of major conventional assets, including ordinary government, investment-grade and high-yield corporate bonds, general stocks, crude oil and gold, from a multiscale perspective. The basic idea is to find out whether green bonds and green equities can be classified as an independent asset class with many common characteristics because of their eco-friendly nature or whether, to the contrary, these green financial products exhibit a greater similarity with their respective non-green equivalents, i.e., regular bonds and general stocks. To this end, a novel approach that integrates wavelet methods and network analysis is applied. 
Network analysis has proven to be a powerful and flexible tool to characterize any complex system as a network of nodes and edges. There is a plethora of works showing that network theory provides an ideal setting to describe the dynamics of information transmission in financial and commodity markets \citep{Man1999,Tse2010,Ji2016,Ji2018}. The networks in this study are constructed based on the wavelet coherence measure of \cite{Gri2004}. Therefore, the proposed wavelet coherence-based network approach allows the level of interdependence between green financial assets and major conventional asset classes to vary across investment horizons.

It is well-known that asset markets in general consist of multiple agents with heterogeneous investment horizons ranging from a few seconds (i.e., intraday traders) to several years (i.e., pension and sovereign wealth funds), who collectively determine aggregate market behavior. Investors with short horizons focus on technical information and are more prone to herding behavior, while investors with longer horizons pay greater attention to fundamental information and the long-term performance of markets \citep{Kri2013}. Accordingly, it seems quite reasonable to think that market participants with diverse horizons will react differently to information flows associated with distinct investment horizons. In this scenario, the level of interdependence between green and and traditional assets can be frequency-dependent and differ greatly among the short-, medium- and long-run.

Understanding the nature of the relationship between green financial instruments and traditional asset classes at diverse time horizons is an issue of undeniable relevance for several economic agents. For investors, who are always searching for alternate sources of return and diversification, it is essential to know the risk-return profile and potential hedging and safe-haven attributes of green investments. This information is particularly valuable in times of market turmoil and financial crisis where there is a great demand for assets that act as safe-haven investments against conventional assets. Moreover, sound knowledge of the interdependence between environmentally friendly financial assets, such as green bonds and green equities, is critical for some investors who want to minimize risk without losing the green properties of their portfolios. This topic is also of great interest to policy makers for several reasons. First, policy makers wish to better understand the intricacies of any emerging market in order to adopt optimal policies. Second, they are committed to develop a sound financial system that makes it possible the raising of the immense amount of capital required for the transition towards a more sustainable economy and, at the same time, exhibits certain resilience to external shocks. 

The main contribution of this paper to the literature resides in that it is, to the best of our knowledge, the first work that investigates the interactions between green financial securities and a number of major conventional asset classes in the time-frequency space using a wavelet-based network approach.
Another important contribution is that it places special emphasis in gaining a profound insight into the linkage between green bonds and green equities, considering not only clean energy stocks as in prior work, but also non-energy green stocks related to green building, green transportation and energy efficiency.

Our empirical results show a close interdependence between green bonds and Treasury and investment-grade corporate bonds, particularly since 2014 coinciding with the start of the exponential growth in green bond issuance, regardless of the investment horizon. Likewise, green stocks are highly tied to general stocks over the entire sample period and irrespective of the time horizon. In contrast, despite their common eco-friendly nature, there is no a strong link between green bonds and green stocks. This means that green bonds and green equities are not perceived by market participants as assets belonging to the same category. In essence, both green bonds and stocks mirror the performance of their respective mainstream equivalents. The low association between green bonds and green equities offers good diversification opportunities to eco-conscious investors without giving up their environmental commitment. Interestingly, a significant connection is found between gold and green bonds, which shows the potential safe-haven properties of green bonds. In addition, there seems to be a weak association between green financial instruments, especially green bonds, and crude oil prices during most of the sample period.

The remainder of the paper proceeds as follows. Section \ref{sec2} contains a brief review of the most relevant literature on the relationship between green financial assets and conventional asset classes. Section \ref{sec3} introduces the wavelet-based network methodology applied in this study, while Section \ref{sec4} describes the data used. Section \ref{sec5} presents and discusses the main empirical results, whereas Section \ref{sec6} concludes the paper.

\section{Literature review}
\label{sec2}

Given the emergent status of the green bond market and that its strongest growth has occurred since 2014 following the publication of the GBP, the academic literature on green bonds is still quite short. Pioneer studies in this area focused on whether green bonds are priced at a premium compared with non-green bonds, providing mixed evidence. For example, \cite{Ehl2017, Bak2018, Feb2018, Gia2019, Zer2019} find a negative premium in green bonds, which implies that investors in green bonds are willing to compensate the environmental benefits of these bonds accepting a lower return. Conversely, according to \cite{Kar2017,Bac2019,Kan2020}, there is a positive yield differential for green bonds. In addition, \cite{Hac2018,Lar2020,Fla2021} find that the green bond premium is essentially zero. A second body of work has demonstrated the benefits of green bond issuance for investors \citep{Bau2019, Tan2020, Wan2020} and issuers \citep{Fla2020, Leb2020}.

The linkage between green bonds and ordinary bonds has been also the subject of several interesting contributions. In a seminal paper, \cite{Pha2016} investigates the volatility behavior of the global green bond market and the global conventional bond market using a multivariate GARCH framework. Her results reveal significant volatility spillovers from the ordinary bond market to the green fixed income market. Subsequently, \cite{Bro2019} utilize a dynamic model averaging method to identify the key determinants of the pattern of correlation between US green and non-green bond markets. They show that the dynamic correlation is affected by a number of factors, such as news-based sentiment towards green bonds, changes in financial market volatility, economic policy uncertainty, crude oil prices and economic activity. In a very recent study, \cite{Nae2021} examine the comparative efficiency of green and conventional bonds employing asymmetric multifractal analysis. Their results show greater inefficiency in the green bond market, especially during downward market trends.

There is also a recent and rapidly growing strand of research that explores the nexus between green bonds and a number of relevant financial and energy asset classes. In this regard, \cite{Reb2018} looks into the co-movement of the global green bond market and a set of mainstream financial and energy markets using bivariate copula models. He concludes that green bonds are tightly connected with Treasury and investment-grade corporate bonds, but tenuously linked to stock and energy commodity markets. Along the same lines, \cite{Reb2020a} analyze the price transmission among green bonds and diverse traditional markets, such as the government, investment-grade and high-yield corporate bond markets as well as stock, exchange rate and energy markets, employing a structural vector autoregressive (VAR) model. The empirical findings show that the green bond market is significantly associated with the government bond and U.S. currency markets and also, but to a smaller degree, to the investment-grade corporate bond market. Specifically, green bonds receive substantial spillovers from these markets, but, however, transmit negligible effects other way around. In contrast, the green bond income market is weakly connected with the high-yield bond, equity and energy markets. Similarly, \cite{Reb2020b} investigate network connectedness across green bonds and several asset classes at different time horizons in the European Union and the U.S. markets. To this end, the spillover index approach of \cite{Die2012,Die2014} is applied over various time scales employing wavelet analysis for decomposing the original series. These authors provide evidence of substantial connectedness across green bonds and government and investment-grade bonds for all horizons in the European Union and the US, even though the relationship is slightly stronger in the short-term. Contrarily, green bonds are weakly linked to high-yield corporate bonds, general equities and energy equities irrespective of the time horizon and the geographical area. It is also shown that oil price uncertainty and, to a minor extent, stock market uncertainty have a significant impact on net spillovers involving green bonds for all the time scales for the European Union and the US.

From an even wider perspective, \cite{Le2021} analyze time and frequency connectedness in volatility across fintech stocks, cryptocurrencies, green bonds and a number of traditional assets, such as general stocks, the US dollar, crude oil, gold and the VIX volatility index, using the Diebold-Yilmaz approach and the connectedness framework of \cite{Bar2018}. Their findings reveal that the transmission of volatility across these markets occurs primarily in the short-term and green bonds appear as net receivers of volatility shocks. Similarly, \cite{Huy2020} examine tail dependence and volatility connectedness across a variety of traditional and novel asset classes, including stocks of artificial intelligence and robotics companies, general stocks, green bonds, gold, crude oil, Bitcoin as well as the VIX index, in the context of the fourth industrial revolution. Copula functions and the connectedness methodology of \cite{Bar2018} are utilized for this purpose. These authors find that portfolios comprised of these asset classes exhibit heavy-tail dependence, implying a high probability of large joint losses in times of financial and economic turbulence. Moreover, transmission of volatility is stronger in the short-term than in the long-term and Bitcoin and gold play a crucial role in hedging.

Analogously, an extensive body of work addresses the interdependence among clean energy stocks, oil prices and other assets, such as technology stocks, general stocks, interest rates, carbon prices and exchange rates, using a variety of econometric techniques \citep{Hen2008, Kum2012, Man2013, Reb2015, Reb2017, Fer2018, Lun2018, Bou2019, Nae2020, Dut2020}. Although there is no general consensus in the literature, most of these studies conclude that stock prices of renewable energy stocks are more strongly related to stock prices of technology companies than to crude oil prices. From a disaggregate perspective, \cite{Pha2019} analyzes whether the linkage between crude oil prices and alternative energy stocks is homogeneous across sub-sectors of the renewable energy equity market. The empirical results show that such relationship varies greatly across different sub-sectors of the alternative energy stock market.

Lastly, another recent strand of literature focuses on the link between green bonds and green equities. From a bivariate perspective, \cite{Ngu2020} investigate the interactions between green bonds and other asset classes, such as general and clean energy stocks, regular bonds and commodities, over time and across investment horizons using the rolling window wavelet correlation approach and the wavelet coherence. They report a significant co-movement between green and ordinary bonds and, although to a lesser extent, between green bonds and alternative energy stocks across practically all time periods and frequencies. There is also a very strong positive correlation between renewable energy stocks and general stocks. However, a weak correlation is observed between green bonds and general stocks and commodities, confirming the safe-haven potential of green bonds. Using copula models and Conditional Value-at-Risk (CoVaR) techniques,\cite{Liu2021} explore the dynamic bivariate dependence structure and risk spillovers between the green fixed income market and a number of global and sectoral renewable energy markets. Positive time-varying average and tail dependence are found between the green bond and the renewable energy stock markets. Furthermore, there is an asymmetric risk spillover between extreme movements in both markets, being the spillover greater for downside risks. Applying the Baruník-Křehlík connectedness framework in a multivariate setting, \cite{Fer2021} analyze the time-frequency connectedness across green bonds and several key mainstream assets, such as Treasury and investment-grade corporate bonds, general and renewable energy stocks, the US dollar exchange rate and crude oil. The results reveal that connectedness primarily occurs at shorter time horizons. Connectedness is particularly strong between green bonds and government and investment-grade corporate bonds, implying that green fixed income securities cannot be regarded as a different asset class. In contrast, there is a quite low connectedness between green fixed income securities and clean energy stocks.
\cite{Ham2020} analyze the time-varying causal relationships between green bonds and a number of financial and environmental variables. These authors find significant causality running from US 10-year Treasury bonds to green bonds, but no causal effects from green bonds to the other variables.

There is also an emerging line of research focused on the interactions between green bonds and traditional asset classes under diverse market circumstances. In this regard, utilizing quantile-based connectedness measures, \cite{Sae2021} analyze return spillovers between green bonds, renewable energy stocks and several dirty energy assets, such as crude oil and stocks of oil and gas companies, under mean and extreme market conditions. They find that spillovers are particularly intense in the lower and upper quantiles. By the same token, \cite{Pha2020} applies the cross-quantilogram approach of \cite{Han2016} to evaluate the dependence between green bond and green stock markets across different quantiles. Her findings indicate that the level of dependence between green bonds and green equities is relatively low under normal market states once controlled for general conditions in the bond, stock and energy markets. However, the dependence rises substantially during extreme market movements.
In another interesting contribution, \cite{Sae2020} examine the hedging ability of green bonds and clean energy stocks against dirty energy assets. They document that renewable energy stocks act as a better hedge than green bonds, particularly in the case of crude oil.

\section{Methodology}
\label{sec3}

In this section, we describe the main features of the novel approach employed in this paper to examine the interdependence between green financial instruments and conventional asset classes. This approach is based on the combination of wavelet methods and network analysis. A basic aspect in constructing networks among assets or markets is the determination of the distance among such assets or markets. Correlation measures are typically used for that purpose and the traditional Pearson correlation coefficient constitutes the most popular indicator. However, the Pearson coefficient is limited by its inability to capture the directionality of the relationships between asset classes or markets and does not take into account that the links can vary across frequencies. So, in an attempt to address these shortcomings, the wavelet coherence measure of \cite{Gri2004} is utilized in this study to calculate the correlations for all the asset pairs at each time point and frequency. Then, the network of green financial products and major traditional assets is constructed based on the estimated wavelet coherence.

\subsection{Wavelet coherence}

Wavelet analysis, a well-established signal processing technique that was born in the mid-1980s \citep{Gou1984,Gros1984}, offers a natural platform to study relationships among time series over time
and across different frequencies simultaneously. In this framework, any signal can be
represented in terms of wavelets. A wavelet is a function $\psi\in L^2(\mathbb{R})$ with zero average (i.e., $\int_{\mathbb{R}}\psi = 0$), normalized ($\|\psi\|=1$) and ``centered'' in the neighborhood of $t=0$ \citep{mallat1998}.  This mother wavelet can be scaled by $s>0$ and translated by $u\in\mathbb{R}$ in order to define the daughter wavelets $\psi _{u,s}(t)$ as:
\begin{equation}
\psi_{u,s}(t) = \frac{1}{\sqrt{s}}\psi \left( \frac{t-u}{s}\right).
\end{equation}
The normalization factor $1/\sqrt{s}$ guarantees unit variance and allows the comparison of wavelet transformations across scales and time.

The Morlet wavelet, given by $\psi(t) = \pi^{-1/4} e^{\mathrm{i}\omega _0 t} e^{-t^2/2}$ \citep{Gou1984} is widely used in economic and financial applications. There, the positive parameter $\omega _0$ denotes the dimensionless frequency of the wavelet. Usually, $\omega _0=6$ has proven an appropriate choice since it provides a good balance between time and frequency localization \citep{Agu2013}.

A major feature of wavelet methods is that they allow the analysis of the relationships among time series over many time scales through cross-wavelet tools such as the wavelet coherence and the wavelet phase difference.  The continous wavelet transform (CWT), $W_x(u,s)$, is obtained by projecting the original time series $x$ onto the mother wavelet as follows:
\begin{equation}
W_x (u,s) = \int _{-\infty}^{+\infty} x(t) \psi_{u,s}^*(t)\,\textrm{d}t ,
\end{equation}
where the superscript $ ^*$ indicates a complex conjugate.

The cross-wavelet transform of two time series x and y is given by:
\begin{equation}
W_{xy} (u,s)=W_x (u,s) W_y^*(u,s),
\end{equation}
where $W_x(u,s)$ and $W_y(u,s)$ are the CWT of $x$ and $y$, respectively.

In turn, the cross-wavelet power spectrum can be computed is the modulus of the cross-wavelet transform, $| W_{xy}\left( u,s\right) |$. It measures the local covariance between the two series for each time and frequency.

The squared wavelet coherence of two time series $x$ and $y$ \citep{Tor1999} is given by:
\begin{equation}
\label{eq:wc}
R^2_{xy} (u,s)=\frac{|S\left( s^{-1} W_{xy} (u,s) \right) |^2}{S\left( s^{-1} |W_x (u,s)|^2 \right) S\left( s^{-1} |W_y (u,s)|^2 \right) },
\end{equation}
being $S$ a smoothing operator both in time and scale, ensuring that the coherence is not constantly equal to one. 

The squared wavelet coherence can be viewed as a measure of local linear correlation both in time and frequency between two signals, analogously to the squared correlation coefficient used in linear regression. In other words, the squared wavelet coherence acts as a correlation coefficient around each moment in time and for each frequency, thus providing a broader picture of the co-movement between two time series than pure time-domain methods. It takes values between 0 and 1, indicating the strength o. 

In order to identify the direction of the relationships between time series, \cite{Tor1999} introduced the \textit{wavelet phase-difference} between two time series $x$ and $y$ at time $u$ and scale $s$, defined by
\begin{equation}
\label{eq:pd}
\phi_{xy} (u,s)=\tan^{-1}\left( \frac{\Im \left( S\left( s^{-1} W_{xy}(u,s) \right) \right) }{\Re \left( S\left( s^{-1} W_{xy}(u,s) \right) \right) }\right) ,
\end{equation}
where $\Re $ and $\Im $ are the real and imaginary parts, respectively. The wavelet phase-difference \eqref{eq:pd} takes values in the interval $\left] -\pi, \pi \right] $ and provides information on lead-lag effects as well as the sign of the association between the two time series. The closer to zero the phase-difference is, the more together both signals move at the specified time and scale. A phase-difference between $-\pi /2$ and $\pi /2$ means that the two time series are in phase (i.e. they move in the same direction), otherwise they are in anti-phase (i.e. they move in the opposite direction). In addition, a positive phase-difference implies that $x$ leads $y$, while a negative phase-difference indicates that $y$ leads $x$.

Based on \eqref{eq:wc} and \eqref{eq:pd}, we can define a pair of ``oriented'' squared wavelet coherences that include the information provided by the phase-difference as:
\begin{equation}
\label{eq:r2xy}
R^2 _{x\rightarrow y}(u,s)=\left\{
\begin{array}{ll}
R^2 _{xy}(u,s) & \textrm{ if }\phi_{xy}(u,s)\geq 0 \\  \\
R^2 _{xy}(u,s) \cdot \left| \cos\left( \phi_{xy}(u,s)\right) \right| & \textrm{ if } \phi_{xy}(u,s)<0,
\end{array}
\right.
\end{equation}
and
\begin{equation}
\label{eq:r2yx}
R^2 _{y\rightarrow x}(u,s)=\left\{
\begin{array}{ll}
R^2 _{xy}(u,s)\cdot \left| \cos\left( \phi_{xy}(u,s)\right) \right| & \textrm{ if }\phi_{xy}(u,s)\geq 0 \\ \\
R^2 _{xy}(u,s) & \textrm{ if } \phi_{xy}(u,s)< 0.
\end{array}
\right.
\end{equation}
Thus, if the phase-difference is positive, then $R^2_{x\rightarrow y}$ maintains the value of $R^2_{xy}$, while $R^2_{y\rightarrow x}$ is penalized, and vice versa. For the particular case of phase-difference equal to $\pi /2$, we have that $R^2_{y\rightarrow x}$ vanishes, whereas if the phase-difference is $-\pi /2$, then $R^2_{x\rightarrow y}$ vanishes.

\subsection{Construction of the network}

A network structure $N$ is defined by a set of nodes (or vertices) $V$ and a set of edges (or links) $E$ joining the nodes, and will be denoted by $N(V,E)$. The universe of asset classes (green and conventional assets) in this study is depicted by nodes. The edges connecting vertices represent the interdependence relationships between each pair of assets, which are measured through the squared wavelet coherence estimated in the previous section. The construction of the network of green financial products and conventional asset classes is based on the value of the wavelet coherence. In particular, this measure of similarity is transformed into a dissimilarity indicator for each time $u$ and scale $s$ simply by calculating the complementary of the squared wavelet coherence as
\begin{equation}
\label{eq:dxy}
D_{xy}(u,s)=1-R^2_{xy}(u,s),
\end{equation}
which obviously also ranges between 0 and 1. 

From \eqref{eq:dxy}, a dissimilarity matrix $D_{N\times N}=\left[ D_{ij}\right]$, containing the dissimilarities between each pair of assets, is defined. Next, a clustering algorithm is applied to this matrix in order to yield clusters of assets that exhibit a synchronized behavior. To this aim, we first determine the optimal number of clusters for each time and scale by means of the Gap statistics \citep{Tib2001}. We have used a Monte Carlo simulation with $50$ time series following a uniform random distribution with the same range than the data at each timestamp. As for the clustering algorithm, we employ the K-medoids algorithm with the PAM (Partitioning around medoids) implementation \citep{Sch2019}. This algorithm has been found to be suitable for time-series clustering in comparison to other more classical methods, such as the K-means clustering \citep{Nie2007,Dur2013,Nak2020}.

Note that this methodology provides a network for each different time $u$ and scale $s$. To construct a single network for a given time interval  $\left[ u_1,u_2\right] $ and a given scale interval $\left[ s_1, s_2\right]$, we have to average the squared wavelet coherence as follows:
\begin{equation}
\label{eq:r2media}
\frac{1}{(u_2-u_1)(k_2-k_1)}\int _{k_1}^{k_2}\int _{u_1}^{u_2} R_{xy}^2 (u,2^k) \, \mathrm{d}u \, \mathrm{d}k ,
\end{equation}
where $k_i=\log_2s_i$ are called \textit{log-scales} \citep{Fer2016}, for $i=1,2$. The fact that the family of dyadic wavelets defines an orthonormal basis of $L^2\left( \mathbb{R}\right)$ \citep{mallat1998} justifies the change to the log-scales, since therefore $k=\log _2s$, and not $s$, is the ``natural'' variable for the scale.

For a given cluster $C$, we can consider that \eqref{eq:r2media} represents the average ``strength'' of the connection between nodes $x$ and $y$, both in $C$, in a given time interval and a given scale interval or frequency band. Furthermore, we can include in \eqref{eq:r2media} the directionality of the connection using the oriented wavelet coherence measures as defined in \eqref{eq:r2xy} and \eqref{eq:r2yx}. Specifically, $R^2_{x\rightarrow y}$ and $R^2_{y\rightarrow x}$ represent the strength of the causal link from $x$ to $y$ and from $y$ to $x$, respectively.

\section{Data description}
\label{sec4}

The dataset used in this research consists of daily closing prices of green bonds and green equities as well as of a number of asset classes that represent well-known investment alternatives. The sample period runs from October 13, 2010 to November 13, 2020, for a total of $2541$ daily observations. The initial date of this sampling period is determined by the availability of data on the green stock market. The variables utilized are described below.  

The S\&P Green Bond index (SPGBI) is used as a proxy for the green bond market \citep{Pha2016, Ngu2020, Huy2020, Le2021}. This index has historical data available since November 2008 and is one of the pioneering indices that track the financial performance of the global green bond market. The SPGBI is a market value-weighted index including only bonds flagged as ``green'' by Climate Bonds Initiative (CBI). These bonds are issued by multilateral, government or corporate entities from any country and denominated in any currency, although the SPGBI is calculated in US dollars.

Various sectoral clean energy indices from the NASDAQ OMX Green Economy family are considered in this study to account for the financial performance of the global green stock market at the sector level. This index family comprises a wide range of indices that capture the behavior of many green equity sub-sectors, thus providing a comprehensive picture of the green stock market. Such indices were originally launched in October 2010 with an initial value of $1000$ and are calculated in US dollars. Most of the proceeds from global green bond issuances between 2010 and 2019 were allocated to Alternative energy (23\% of total issuance volume in US dollars), Energy efficiency (19\%), Clean transportation (14\%), Sustainable water management (12\%), Green buildings (9\%) and Pollution prevention (7\%) projects \citep{Ire2020}. In order to maintain the consistency between the data of green bonds and green stocks, only the equity indices from the NASDAQ OMX Green Economy family corresponding to the sectors with larger relative importance in terms of green bond issuance are considered in this paper\footnote{The only exception is the Sustainable water management sector. This sector has not been considered in the empirical analysis because the time series of the corresponding green equity index is not available since the beginning of the sample period.}. Thus, we have utilized the Solar Energy (GRNSOLAR), Wind Energy (GRNWIND), Fuel Cell (GRNFUEL), Energy Efficiency (GRNENEF), Clean transportation (GRNTRN), Green buildings (GRNGB) and Pollution prevention (GRNPOL) equity indices. The GRNSOLAR and GRNWIND indices are designed to track the performance of companies that produce energy through solar and wind power, respectively, while the GRNFUEL index is comprised of firms that generate energy through fuel cells. The GRNENEF tracks companies that dramatically increase energy efficiency and the GRNTRN is composed of companies focused on efficiency gains and pollution reduction in automobiles, trains and other means of transport. Lastly, the GRNGB index reflects the performance of companies participating in advanced designs for retrofits and new buildings that lead to dramatic efficiency gains in energy and water consumption, while the GRNPOL index tracks firms producing goods and services that reduce pollution generated by traditional industries. It is worth mentioning that this selection includes both renewable energy and non-energy sectors of the green equity market.

As for the conventional asset classes, the global Treasury, investment-grade and high-yield corporate bond markets are proxied by the \textit{Bloomberg Barclays Global Treasury total return index}, the \textit{Bloomberg Barclays Aggregate Corporate bond index} and the \textit{Bloomberg Barclays Global High Yield Corporate total return index}, respectively. In turn, the \textit{MSCI World index} is utilized to characterize the performance of the overall stock market. Furthermore, the \textit{Brent oil price} is employed as a proxy for the international crude oil market. Oil price fluctuations can affect the dynamics of green bonds and equities as lower crude oil prices reduce the incentives to substitute oil by clean energy, thus adversely impacting the performance of green financial assets. In addition, the \textit{S\&P GSCI gold index} is utilized to measure the spot price of gold. Gold has traditionally acted as a safe haven for investors in stocks, bonds and currencies during turbulent periods. Hence, it is interesting to examine whether there is a certain connection between green financial assets and gold, particularly in times of market stress or turmoil.

All data are collected from Thomson Reuters Datastream. Following the usual practice, daily returns calculated for each asset class as the logarithmic difference of the closing price between two consecutive observations are used in the empirical analysis.

Table~\ref{tab:descrstat} reports the main descriptive statistics of the daily returns of all series over the whole sample period. The mean daily returns are close to zero for all the series and quite smaller than their respective standard deviations, which suggests relatively high variability of all asset classes under scrutiny. Crude oil returns have a very high standard deviation, only surpassed by the Fuel Cell equity sub-sector. The high volatility of the oil market can be attributed to the tremendous swings in the price of oil since the beginning of the 21st century. As expected, the green bond index has smaller return and standard deviation than the green equity sector indices and the overall stock market. The skewness coefficient is negative for the vast majority of return series, indicating that the distributions are left-skewed. In turn, the kurtosis values are higher than 3 for all return series, implying a leptokurtic distribution with fatter tails than those of a normally distributed series. In addition, the Jarque-Bera (JB) test statistics reject the null hypothesis of normality in all cases at the 1\% level, corroborating the departure from normality. The standard ADF unit root and KPSS stationarity tests are also applied to check the unit root properties of the series. Their results indicate that all return series are stationary processes, i.e., $I(0)$, at the 1\% significance level.

\begin{table}[]
\centering
{\footnotesize
\caption{Descriptive statistics of the data. This table presents the descriptive statistics and unit root test results of daily returns of green bonds and green stocks and the group of conventional asset classes considered for the whole sample period. JB denotes the J-B test statistic for the null of normality. In turn, ADF and KPSS represent the statistics of the Augmented Dickey-Fuller (ADF) unit root test and Kwiatkowski-Phillips-Schmidt-Shin (KPSS) stationarity test, respectively. As usual, $^{***}$ indicates statistical significance at the 1\% level.}\label{tab:descrstat}
\begin{tabular}{lrrrrrrr}
  \toprule
 & Mean & Std & Skn & kurt & JB & KPSS & ADF \\ 
  \midrule
GOLD & 0.055 & 4.517 & -0.622 & 6.958 & 5299.012$^{***}$ & 0.085 & -16.847$^{***}$ \\ 
GRNENEF & 0.128 & 5.668 & -0.467 & 9.214 & 9095.839$^{***}$ & 0.350 & -16.589$^{***}$ \\ 
GRNFUEL & 0.218 & 13.515 & 0.285 & 5.913 & 3743.917$^{***}$ & 0.133 & -17.326$^{***}$ \\ 
GRNGB & 0.030 & 5.779 & -1.805 & 27.108 & 79288.121$^{***}$ & 0.186 & -16.799$^{***}$ \\ 
GRNPOL & 0.129 & 5.027 & -0.823 & 9.315 & 9489.666$^{***}$ & 0.445 & -17.328$^{***}$ \\ 
GRNSOLAR & 0.150 & 8.605 & -0.593 & 7.716 & 6464.576$^{***}$ & 0.196 & -16.030$^{***}$ \\ 
GRNTRN & 0.242 & 5.709 & -0.853 & 15.598 & 26107.321$^{***}$ & 0.546 & -16.003$^{***}$ \\ 
GRNWIND & 0.234 & 6.890 & -0.309 & 3.202 & 1128.656$^{***}$ & 0.384 & -15.558$^{***}$ \\ 
HYBOND & 0.100 & 1.381 & -2.273 & 34.402 & 127657.710$^{***}$ & 1.638 & -13.941$^{***}$ \\ 
IGBOND & 0.066 & 1.279 & -1.460 & 15.518 & 26437.824$^{***}$ & 0.980 & -16.719$^{***}$ \\ 
MSCI & 0.124 & 4.186 & -1.164 & 17.186 & 31889.063$^{***}$ & 0.601 & -16.440$^{***}$ \\ 
OIL & -0.121 & 12.647 & -3.291 & 122.720 & 1601031.621$^{***}$ & 0.092 & -14.438$^{***}$ \\ 
SPGBI & -0.012 & 1.752 & -0.356 & 5.865 & 3702.540$^{***}$ & 0.337 & -18.027$^{***}$ \\ 
TREAS & 0.028 & 1.584 & -0.283 & 3.434 & 1285.460$^{***}$ & 0.090 & -18.277$^{***}$ \\ 
   \bottomrule
\end{tabular}
}
\end{table}

Table~\ref{tab:correlation} shows the Pearson correlation coefficients between each pair of assets over the full sample period. Not surprisingly, a very high positive correlation is found between green bonds and conventional Treasury ($0.72$) and investment-grade corporate bonds ($0.70$), which indicates an important linkage between green bonds and conventional high credit-quality bonds. Interestingly, the unconditional Pearson correlations between green bonds and green stocks are rather heterogeneous and not very high, depending largely on the particular green equity sub-sector considered. More precisely, these correlations range from $0.13$ for the Fuel Cell sub-sector to $0.41$ for the Pollution Prevention sub-sector. This relatively low correlation suggests that overall green bonds and green equities do not seem to be driven by the same factors. In fact, the green equity sector indices tend to be more significantly correlated among themselves and the general stock market than with the green bond index. Furthermore, a very low positive correlation ($0.10$) is observed between green bonds and crude oil prices. Likewise, green stocks exhibit a relatively low positive correlation (ranging from $0.19$ to $0.30$) with crude oil, though more pronounced than in the case of green bonds. Lastly, it should be mentioned that green stocks are in general weakly correlated with Treasury and investment-grade corporate bonds and with gold.

\begin{landscape}

\begin{table}[]
\centering
{\scriptsize
\caption{Correlation matrix. This table reports the values of the Pearson linear correlation coefficient between all pairs of variables over the full sample period.}\label{tab:correlation}
\begin{tabular}
{lrrrrrrrrrrrrrr}
  \toprule
 & GOLD & GRNENEF & GRNFUEL & GRNGB & GRNPOL & GRNSOLAR & GRNTRN & GRNWIND & HYBOND & IGBOND & MSCI & OIL & SPGBI & TREAS \\ 
  \midrule
GOLD & 1.00 &  &  &  &  &  &  &  &  &  &  &  &  &  \\ 
  GRNENEF & 0.06 & 1.00 &  &  &  &  &  &  &  &  &  &  &  &  \\ 
  GRNFUEL & 0.05 & 0.42 & 1.00 &  &  &  &  &  &  &  &  &  &  &  \\ 
  GRNGB & 0.08 & 0.77 & 0.39 & 1.00 &  &  &  &  &  &  &  &  &  &  \\ 
  GRNPOL & 0.11 & 0.77 & 0.35 & 0.70 & 1.00 &  &  &  &  &  &  &  &  &  \\ 
  GRNSOLAR & 0.08 & 0.70 & 0.44 & 0.60 & 0.55 & 1.00 &  &  &  &  &  &  &  &  \\ 
  GRNTRN & 0.10 & 0.72 & 0.44 & 0.67 & 0.61 & 0.64 & 1.00 &  &  &  &  &  &  &  \\ 
  GRNWIND & 0.07 & 0.51 & 0.28 & 0.48 & 0.56 & 0.47 & 0.43 & 1.00 &  &  &  &  &  &  \\ 
  HYBOND & 0.17 & 0.51 & 0.29 & 0.61 & 0.56 & 0.43 & 0.50 & 0.51 & 1.00 &  &  &  &  &  \\ 
  IGBOND & 0.41 & 0.04 & 0.03 & 0.21 & 0.17 & 0.05 & 0.11 & 0.20 & 0.53 & 1.00 &  &  &  &  \\ 
  MSCI & 0.10 & 0.92 & 0.45 & 0.84 & 0.80 & 0.73 & 0.78 & 0.57 & 0.64 & 0.13 & 1.00 &  &  &  \\ 
  OIL & 0.11 & 0.30 & 0.21 & 0.27 & 0.29 & 0.26 & 0.30 & 0.19 & 0.30 & 0.03 & 0.34 & 1.00 &  &  \\ 
  SPGBI & 0.43 & 0.31 & 0.13 & 0.36 & 0.41 & 0.22 & 0.29 & 0.34 & 0.53 & 0.70 & 0.37 & 0.10 & 1.00 &  \\ 
  TREAS & 0.46 & -0.04 & -0.02 & 0.06 & 0.07 & -0.04 & 0.01 & 0.07 & 0.26 & 0.81 & 0.01 & -0.04 & 0.72 & 1.00 \\ 
   \bottomrule
\end{tabular}
}
\end{table}
\end{landscape}

\section{Empirical results}
\label{sec5}

In this section, we first present and discuss the network graphs of green financial instruments and conventional asset classes based on the squared wavelet coherence for the entire sample period. Then, we check the robustness of the results by examining the network graphs for several sub-periods of heightened economic and financial uncertainty. Three different frequency bands are considered to assess the degree of interdependence between green investments and mainstream assets across various investment horizons. Specifically, the first frequency band refers to a horizon between $2$ and $5$ days, i.e., one business week, and, obviously, represents the short-term. The second spectral band covers from $5$ days until $22$ days, i.e., approximately one business month, and is considered adequate to capture the behavior of financial and commodity markets in the medium-term. Lastly, the third frequency band represents an investment horizon of more than $22$ days and is associated with the long-term.

\subsection{Full sample analysis}

\begin{figure}[]
\center
\includegraphics[width=10.5 cm]{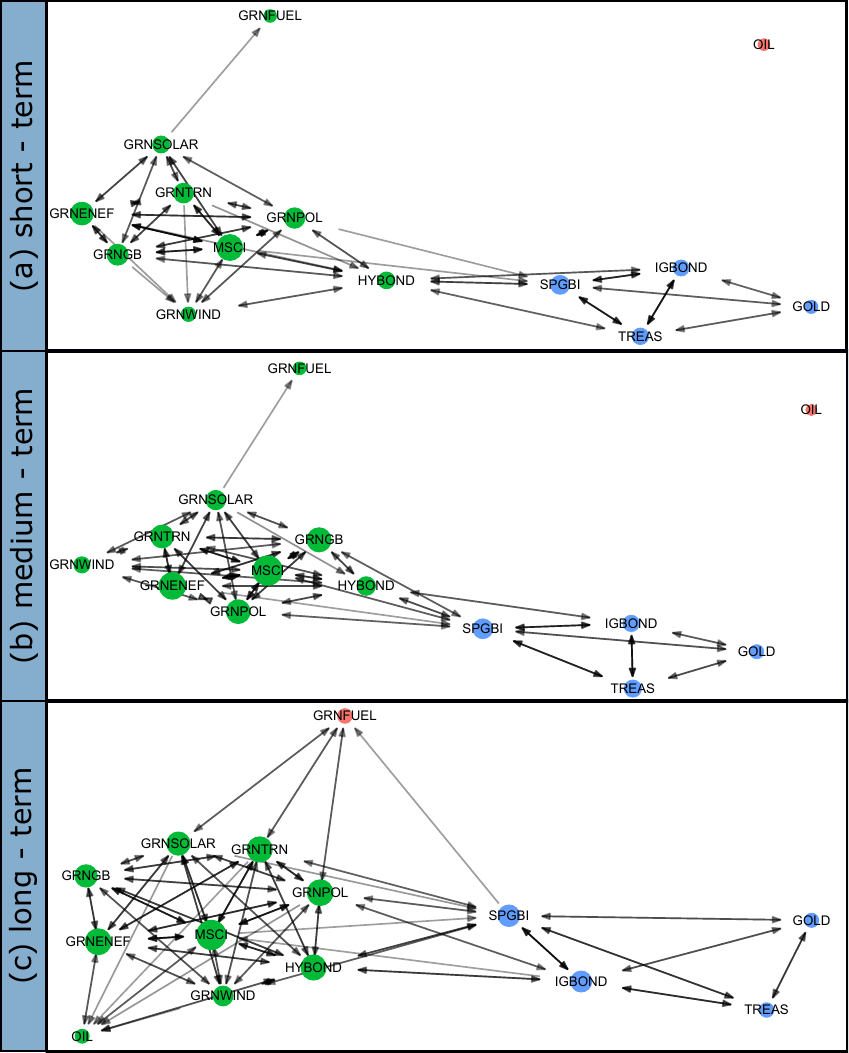}
\caption{Wavelet coherence-based networks over the full sample for the three frequency bands considered. (\textbf{a}) Short-term, between $2$ and $5$ days. (\textbf{b}) Medium-term, between $5$ and $22$ days. (\textbf{c}) Long-term, more than $22$ days. Source: own elaboration using the R statistical software \citep{Rbase} together with packages wavScalogram \citep{wavscalogram}, igraph \citep{igraph}, and ggraph \citep{ggraph}.\label{fig:global}}
\end{figure}

Figure~\ref{fig:global} depicts the squared wavelet coherence-based networks of green financial products and traditional asset classes over the full sample for the three frequency bands considered. Specifically, plots (a), (b) and (c) refer to the short-term, medium-term and long-term, respectively. The networks are comprised of asset classes as nodes and links based on the mean of the squared wavelet coherence as directed edges.
The optimal number of clusters is chosen applying the Gap statistic \citep{Tib2001} to the output of the K-medoids clustering algorithm \citep{Sch2019}, using as a dissimilarity measure the expression in \eqref{eq:dxy}, but replacing the squared wavelet coherence by its average value as calculated in \eqref{eq:r2media}. It is worth mentioning that in this approach, the goodness of a clustering assignation is related to the average intra-cluster distance. The optimal number of clusters is determined by comparing the results obtained with the original data with those resulting from a Monte Carlo simulation of a given number of random signals ($50$ in this case) generated uniformly over the range of each feature (each time observation). In particular, the optimal number of clusters is the minimum number of clusters for which the results of the Monte Carlo simulation are undistinguishable from the results obtained with the real data (for more details see, for instance, \cite{Tib2001}).

The black density of the edges connecting nodes is proportional to the value of the mean of the squared wavelet coherence between them given by \eqref{eq:r2media}. Moreover, the directional information provided by the wavelet phase-difference, as reflected by the arrowheads, is included through the mean of the oriented squared wavelet coherences given by \eqref{eq:r2xy} and \eqref{eq:r2yx}. It is important to remark that we have not used these oriented versions of the squared wavelet coherence in the determination of the number of clusters because the goal was that nodes in the same cluster have a high coherence between them, regardless of the phase-difference. Note that if we had used the oriented versions in \eqref{eq:r2xy} and \eqref{eq:r2yx}, then the edges with strong causality (phase-difference near $\pi/2$ or $-\pi/2$) would be discriminated against edges with weak causality (phase-difference near $0$ or $\pm \pi$) since if there is a strong causality, then the squared wavelet coherence in one direction is severely penalized.

All nodes in the same cluster have the same color. The size of each node is proportional to the number and strength of the connections from the node to all the other nodes in the system. Specifically, the size of a node $x$ is computed as the mean of the oriented squared wavelet coherences $R^2_{x\rightarrow y}$ in \eqref{eq:r2xy} for $y$ being any other node.
Only the edges and arrowheads representing squared wavelet coherences above $0.38$ are graphically displayed. We have chosen this threshold because, according to a Monte Carlo simulation of $100$ repetitions, the mean of the squared wavelet coherence between pairs of Gaussian noise series (with the same variance as the original series) is not greater than $0.38$ with a $95\%$ of confidence in all the frequency bands.

Three quite homogeneous clusters are apparent in Figure~\ref{fig:global} regardless of the investment horizon. The first cluster is mainly composed of the different types of green stocks considered as well as general stocks and high-yield corporate bonds. The second cluster comprises green bonds, Treasury bonds, investment-grade corporate bonds and gold and, because of its composition, it can be categorized as a cluster of low-risk assets. The third cluster is smaller and has a less consistent structure. It is formed only by crude oil in the shorter and medium horizons and fuel cell stocks in longer horizons. Besides belonging to the same cluster, there are strong links between green bonds and Treasury and investment-grade corporate bonds for all investment horizons. A convincing explanation for this result is based on the fact that green and ordinary fixed income securities share many characteristics in terms of credit rating, issuer, currency, coupon rates and maturity. Similar strong relationships between green bonds and government and high-credit rating corporate bonds have been documented in \cite{Reb2018, Ngu2020, Pha2020, Reb2020a, Reb2020b, Fer2021}. Equally, a strong nexus is found between green stocks and general stocks irrespective of the horizon, with the sole exception of fuel cell stocks\footnote{In this respect, it should be noted that \cite{Pha2019} also finds that fuel cell stocks exhibit a quite decoupled behavior from the other clean energy stocks and the crude oil market.}. This substantial connection was totally expected since green equities are an increasingly important part of the overall equity market and is in line with the evidence reported in \cite{Ngu2020, Pha2020, Fer2021}. As can be seen, there exists a remarkable association between green financial products and traditional assets for all horizons, although the intensity of connections seems to increase slightly as the investment horizon becomes longer. This finding is not surprising considering that the relationships among markets at shorter horizons tend to be contaminated by extraordinary transitory events, such as changes in market sentiment, financial panic or other psychological factors.

In contrast, it is worth noting that green bonds and green stocks are part of different clusters for any investment horizon, so that they are not closely tied. This evidence suggests that, despite their common environmentally friendly nature, green bonds and green equities are not perceived as a separate asset class. In fact, green bonds are, above all, a fixed income security and green stocks a type of stock. This means that green bonds are directed at investors with a more conservative profile and greater risk aversion than typical investors in green equities. In this regard, \cite{Ngu2020, Pha2020, Fer2021, Sae2021} also fail to identify an important association between green bonds and clean energy stocks.

Another relevant result is that high-yield corporate bonds, in spite of their fixed income nature, are not closely interlinked with green bonds. Actually, high-yield bonds are in the same cluster than general stocks and green stocks regardless of the investment horizon considered. Therefore, high-yield corporate bonds are more similar to stocks than to Treasury and high-credit quality corporate bonds. This result can be attributed to the fact that high-yield bonds and stocks both have substantial exposure to default risk. The limited connection between green bonds and high-yield bonds is fully consistent with the findings of \cite{Reb2018, Reb2020a, Reb2020b, Fer2021}.

Interestingly, a close linkage is observed between gold and green bonds. In particular, gold is in the same cluster than green bonds for all horizons. The inclusion of green bonds within the cluster of low-risk assets made up of gold, Treasuries and investment-grade corporate bonds indicates that green bonds can play a significant role as diversifiers and safe-haven assets against some riskier investments such as equities, mainly during turbulent times. In the same spirit, \cite{Ari2021} highlights the potential role of green bonds as an effective hedging and safe-haven instrument for traditional assets during the COVID-19 pandemic. In addition, there appears to be a decoupling of crude oil prices from all other asset classes, especially in the short- and medium-term. This finding implies that oil price developments have not played a key role in the impressive expansion of green bond and stock markets during the last few years. Instead, this tremendous growth seems to be more related to the rising environmental awareness of investors, the improved transparency of green financial markets and the good market performance of green bonds and equities.

\subsection{Dynamic analysis}

\begin{figure}[]
\center
\includegraphics[width=10.5 cm]{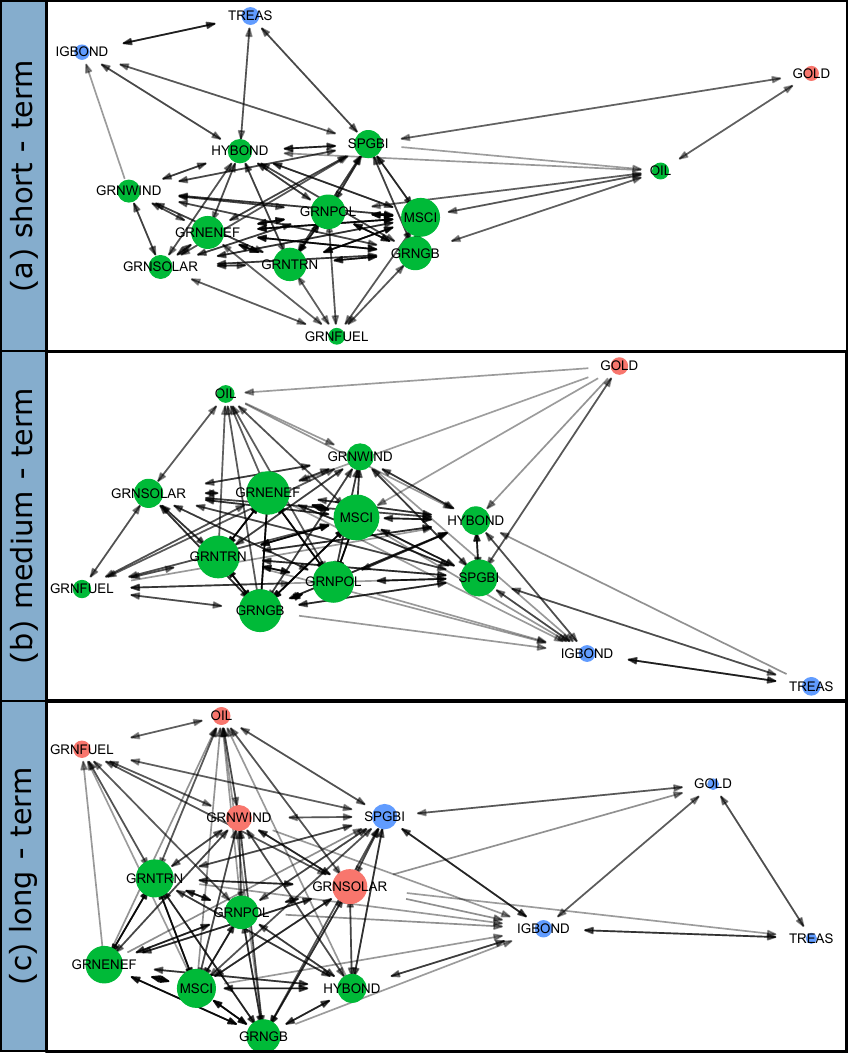}
\caption{Wavelet coherence-based networks during the European sovereign debt crisis sub-period (from October 13, 2010 to July 31, 2012) for the three frequency bands considered. (\textbf{a}) Short-term, between $2$ and $5$ days. (\textbf{b}) Medium-term, between $5$ and $22$ days. (\textbf{c}) Long-term, more than $22$ days. Source: own elaboration using the R statistical software \citep{Rbase} together with packages wavScalogram \citep{wavscalogram}, igraph \citep{igraph}, and ggraph \citep{ggraph}.\label{fig:t1}}
\end{figure}

\begin{figure}[]
\center
\includegraphics[width=10.5 cm]{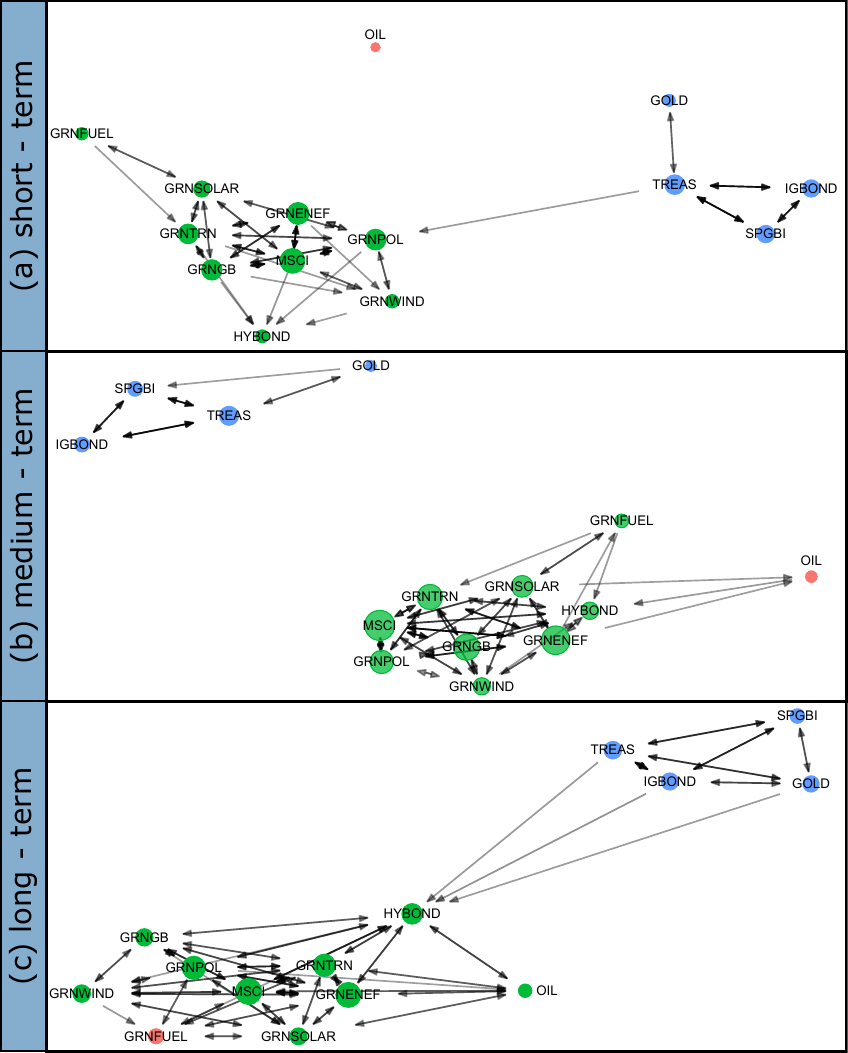}
\caption{Wavelet coherence-based networks during the oil price collapse sub-period (from June 20, 2014 to February 28, 2016) for the three frequency bands considered. (\textbf{a}) Short-term, between $2$ and $5$ days. (\textbf{b}) Medium-term, between $5$ and $22$ days. (\textbf{c}) Long-term, more than $22$ days. Source: own elaboration using the R statistical software \citep{Rbase} together with packages wavScalogram \citep{wavscalogram}, igraph \citep{igraph}, and ggraph \citep{ggraph}.\label{fig:t2}}
\end{figure}

\begin{figure}[]
\center
\includegraphics[width=10.5 cm]{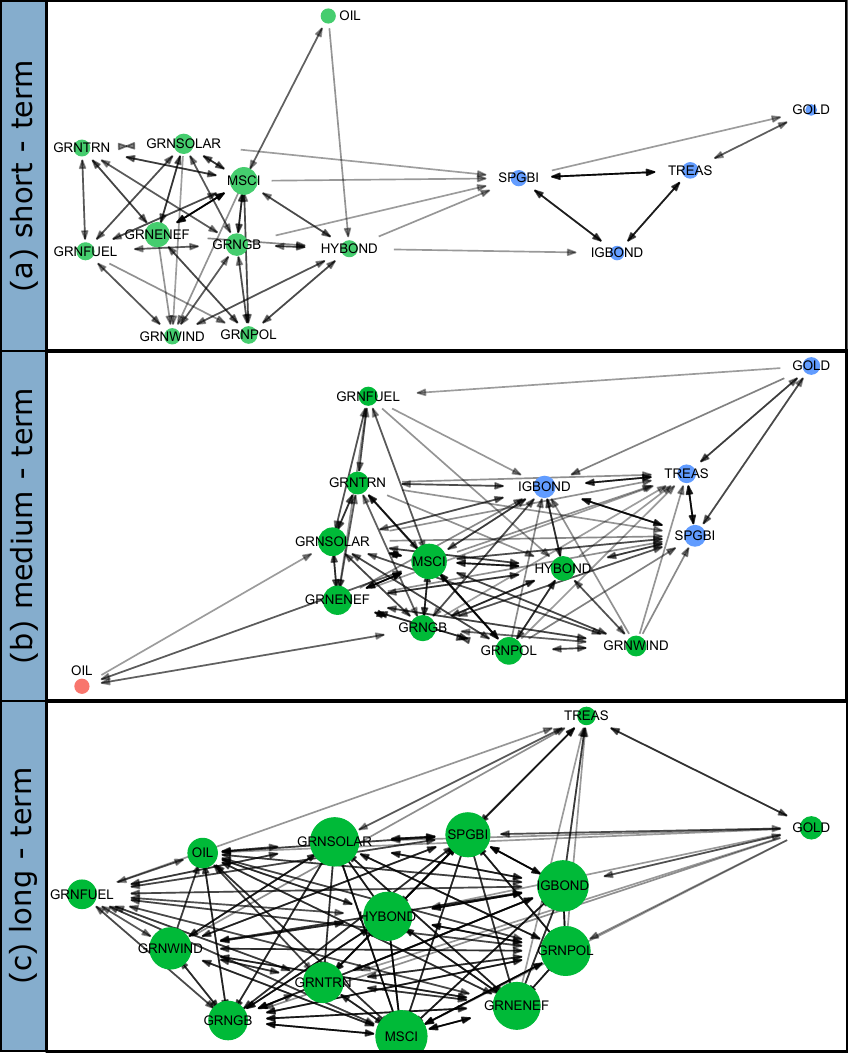}
\caption{Wavelet coherence-based networks during the COVID-19 sub-period (from January 1, 2020 to November 13, 2020) for the three frequency bands considered. (\textbf{a}) Short-term, between $2$ and $5$ days. (\textbf{b}) Medium-term, between $5$ and $22$ days. (\textbf{c}) Long-term, more than $22$ days. Source: own elaboration using the R statistical software \citep{Rbase} together with packages wavScalogram \citep{wavscalogram}, igraph \citep{igraph}, and ggraph \citep{ggraph}.\label{fig:t3}}
\end{figure}

In order to determine whether the level of interdependence between green financial products and traditional asset classes varies significantly during turbulent times, network analysis is carried out separately for various sub-periods associated with some major financial, economic and even health-related events. More precisely, the full sampling period is divided into three sub-periods, namely, the European sovereign debt crisis (from October 13, 2010 to July 31, 2012), the oil price collapse (from June 20, 2014 to February 28, 2016) and the COVID-19 pandemic in 2020 (from January 1, 2020 to November 13, 2020). The Eurozone debt crisis was a period of huge financial turmoil triggered by concerns about the high levels of public debt in several Southern European countries (Greece, Ireland, Italy, Portugal and Spain) and even raised fears about the survival of the Euro. The oil price crash between mid-2014 and early 2006 witnessed a $70\%$ drop in crude oil prices due to a combination of demand and supply factors and also generated considerable turbulence in financial and commodity markets. Lastly, the outbreak of the COVID-19 pandemic caused global financial markets to fall dramatically in March 2020 and to bounce back later as well as a massive decline in economic activity. Networks based on the squared wavelet coherence are built for each sub-period using the methodology described in the previous sub-section. Figure~\ref{fig:t1} displays the squared wavelet coherence-based networks during the European sovereign debt crisis sub-period for the three frequency bands, while Figures~\ref{fig:t2} and \ref{fig:t3} refer to the oil price collapse and COVID-19 sub-periods, respectively.

The results for sub-periods are largely consistent with those of the full sample analysis in terms of the composition of clusters and the most important connections. As a matter of fact, green bonds, regular bonds and gold, on the one hand, and green stocks, general stocks and high-yield corporate bonds, on the other hand, continue to be the two primary clusters and the strongest links are found within these two clusters. However, it is worth mentioning that, during the sub-period associated with the Eurozone sovereign debt crisis, green bonds are more tightly connected to green stocks, mainly in the short- and medium-term, than to government and investment-grade corporate bonds. One possible explanation for this result is that during the European debt crisis the green bond market was still in a very early stage and the perception that green bonds were closer to other green assets, such as green equities, to the detriment of regular bonds was predominant.

Importantly, the links between green and conventional asset classes tend to be more intense during turmoil periods than over the entire sample period regardless of the investment horizon. This is especially visible during the Eurozone sovereign debt crisis and the COVID-19 crisis, in accordance with the commonly held view that cross-market linkages tend to be particularly important in times of market turbulence.
In the same vein, \cite{Bou2021} discover that the structure of connectedness across financial and commodity markets changed drastically around the COVID-19 outbreak.

The sub-period results confirm the lack of a close association between green bonds and green stocks since the green bond market reached certain levels of size and maturity following the publication of the GBP in 2014. The only exception is found in the long-term during the COVID-19 crisis in which a single cluster is identified. However, this evidence should be interpreted with caution because of the border effects and the limited number of available observations at longer horizons owing to the short duration of the COVID-19 sub-period. The sub-period analysis corroborates the existence of a solid link between green bonds and gold, mainly during the oil price collapse and the COVID-19 pandemic, supporting the hedging and safe-haven properties of green bonds during episodes of market turmoil. Additionally, the results for sub-periods show the relative isolation of crude oil from the remaining assets irrespective of the investment horizon in turbulent times, particularly during the oil price crash between mid-2014 and early 2016 and the COVID-19 crisis.

\section{Conclusions}
\label{sec6}

Mounting concerns about climate change and its harmful consequences on human life and biodiversity have resulted in an increased interest from investors and policy makers regarding sustainable finance. 
This paper explores the interdependence between the two most important green financial instruments, namely green bonds and green stocks, and a group of major mainstream asset classes, such as Treasury, investment-grade and high-yield corporate bonds, general equities, crude oil and gold, over various investment horizons. A novel wavelet-based network approach is employed for this purpose. 

Our empirical results clearly show a high level of interdependence between green financial instruments and traditional asset classes regardless of the investment horizon considered. However, the links seem to be slightly stronger in the long-term, consistent with the idea that the relationships between asset markets at longer horizons are less disturbed by ephemeral phenomena such as shifts in investors' market sentiment or psychological factors than shorter-term relationships. A very close association is observed between green bonds and Treasury and investment-grade corporate bonds, mainly from the takeoff of the green bond market in 2014, in the short-, medium- and long-term. This finding can be attributed to the numerous similarities between green bonds and ordinary bonds as far as issuers, currency, credit rating, maturity, coupon rates, etc. Analogously, a strong connection is detected between green stocks and general stocks over the whole sample period irrespective of the time horizon. This is not surprising either taking into consideration that green equities are a significant part of the overall equity market. In contrast, the environmental-friendly nature shared by green bonds and stocks does not translate into a tight linkage between both types of green instruments consistently throughout the sample period and for any horizon. Therefore, green bonds and green equities cannot be seen as an independent asset class, but rather they are two separate assets with distinct risk-return profile and aimed at a different type of investor. Actually, our findings support the idea that the intrinsic character of green bonds and stocks prevails over the commitment toward the environment common to both green assets.

It is also worth highlighting the solid link between green bonds and gold in the short-, medium- and long-term. This evidence is of great interest as it shows that, similar to gold, green bonds can be used as effective hedgers and safe-haven assets, especially during periods of market turmoil. On the contrary, there are no strong ties between green bonds and high-yield corporate bonds, except during the initial stage of the green bond market in the early 2010s. The results also reveal that crude oil is relatively isolated from the remaining asset classes during most of the sample period, particularly in the medium- and long-term. This implies that the behavior of green bond and stock markets over the last few years has been virtually independent of the evolution of oil prices. Furthermore, interdependence between green financial products and traditional asset classes escalates considerably during periods of heightened uncertainty, such as the European sovereign debt crisis and the COVID-19 pandemic.

The empirical evidence reported in this paper has several important implications for economic agents. First, investors should be fully aware that green bonds do not offer significant diversification benefits in portfolios primarily made up of Treasury and investment-grade corporate bonds. Likewise, green stocks also do not work as a good diversifier against general stocks. However, the weak association between green bonds and green stocks implies that climate-conscious investors can benefit from interesting diversification opportunities by including both types of green financial product in their portfolios without compromising their environmental values. In addition, the hedging and safe-haven potential of green bonds against conventional assets can be exploited by investors in times of turbulent markets. 

Second, policy makers, who are interested in the development of a powerful financial sustainable system that facilitates the attainment of governments' ambitious climate-change goals, should assume the importance of policies promoting the expansion of the green bond and green stock markets. Larger, more mature and efficient green bond and equity markets will boost confidence of investors and issuers in these markets and will help reduce their vulnerability during episodes of financial turmoil, helping not to interrupt the mobilization of capital to sustainable investment. In this regard, the lack of a close link between green bonds and green stocks is positive from the standpoint of policy makers as it translates into lower fragility of the green financial system to external shocks. In any case, policy actions supporting the development of green financial instruments must be viewed as a complement, not a substitute, for energy policies developed in the real economy in an attempt to accelerate the transition towards a greener economy.

Finally, this study opens the door to new avenues for future research.
For example, given the apparent hedging and safe-haven properties of green bonds, a possible
line of future research could be aimed at delving deeper into this role through the use of
methods specifically designed for the analysis and identification of asset classes exhibiting
safe-haven features. Another pathway would be to evaluate the diversification properties of
green financial products from a perspective of portfolio and risk management analysis by
constructing portfolios comprised of green and conventional asset classes.

\bibliography{references.bib}
\bibliographystyle{abbrvnat}

\end{document}